\documentclass[12pt]{amsart}
\usepackage{amssymb}

\newtheorem{theorem}{Theorem}
\newtheorem{lemma}[theorem]{Lemma}
\newtheorem{prop}[theorem]{Proposition}
\newtheorem{cor}[theorem]{Corollary} 
\theoremstyle{definition}
\newtheorem{example}[theorem]{Example}

\DeclareMathOperator\cn{\mathfrak c}
\begin{document}

\title{The content of a Gaussian polynomial is invertible}

\author{K. Alan Loper}
\address{Department of Mathematics, Ohio State University-Newark,
Newark, Ohio 43055}
\email{lopera@math.ohio-state.edu}

\author{Moshe Roitman}
\address{Department of Mathematics, University of Haifa,
Haifa 31905, Israel}
\email{mroitman@math.haifa.ac.il}
\thanks{M. Roitman thanks the Mathematics Department of 
the Ohio State University for its hospitality}

\subjclass[2000]{13 B25}
\keywords{content, Gaussian polynomial, invertible ideal, locally principal, prestable ideal}

\begin{abstract}
Let $R$ be an integral domain and let $f(X)$ be a nonzero polynomial 
in $R[X]$.  The content of $f$ is the ideal $\mathfrak c(f)$ generated by the 
coefficients of $f$.  The polynomial $f(X)$ is called 
Gaussian if $\mathfrak c(fg) = \mathfrak c(f)\mathfrak c(g)$ for all
$g(X) \in R[X]$.  It is well known that if $\mathfrak c(f)$ is an invertible ideal, then $f$ is Gaussian.  In this note we prove the converse.  
\end{abstract}

\maketitle

Let $R$ be a ring, that is, a commutative ring with unity. Let $A$ be a ring extension of $R$, and let $f(X) \in A[X]$ be a polynomial:
\[ f(X) = a_{n}X^{n} + \ldots + a_{1}X + a_{0}. \]
The {\em content} ideal of $f$, designated by $\mathfrak c(f)=\mathfrak c_R(f)$,
 is the $R$-submodule of $A$ generated  by the coefficients of $f$ in $A$.  It is easy to see that if $f,g \in R[X]$ (and $A=R$), then 
 $\mathfrak c(fg) \subseteq \mathfrak c(f)\mathfrak c(g)$.
The inclusion in this statement is generally proper.  
The polynomial $f(X) \in R[X]$ is said to be 
{\em Gaussian} if $\mathfrak c(fg) = \mathfrak c(f)\mathfrak c(g)$ 
holds for all  $g(X) \in R[X]$. It is well known that
if $\mathfrak c(f)$ is an invertible ideal of $R$, then $f$ is Gaussian.
More generally, if $R$ is any ring and $\mathfrak c(f)$ is locally principal, then $f$ is Gaussian (see, e.g., \cite[Theorem 1.1]{AK}).
 Recall that a nonzero ideal $I$ of an integral domain $R$ is invertible iff it is locally principal, that is, iff $IR_M$ is a principal ideal for each maximal ideal $M$ of $R$. 

For general background see \cite{G}.

It has been conjectured that the converse is true if $R$ is an integral domain (see  
\cite{{AK}, {GV}}), that is, a Gaussian polynomial over an integral domain has an invertible content. This question is included in the Ph.D. thesis of Kaplansky's student H. T. Tang. 
Significant progress has been 
made on this conjecture in two recent papers \cite{{GV},{HH}}: in \cite{GV}, Glaz and Vasconcelos prove the conjecture for $R$ integrally closed with some additional assumptions (including the Noetherian case). The general Noetherian case is settled in \cite{HH}. As explained below, the conjecture has a local character; thus Heinzer and Huneke prove the conjecture for any locally Noetherian integral domain (this result is obtained as a particular case of a more general theorem).

For connections with the Dedekind-Mertens Lemma see \cite{HH2}. Note that the Dedekind-Mertens Lemma implies that $\sqrt{\mathfrak c (fg)}=\sqrt{\mathfrak c( f)
\mathfrak c(g)}$ for any polynomials $f(X),g(X)$ over an arbitrary ring $R$.

The purpose of this note is to prove the conjecture for 
all integral domains.  
The Gaussian property of a polynomial $f(X)\in R[X]$ is local, that is, $f$ is Gaussian iff the image of $f$ in $R_M[X]$ is Gaussian for each maximal ideal $M$ of $R$. Thus to prove the conjecture we may assume that $R$ is quasilocal (cf. \cite{GV} and \cite{HH}). Moreover, this allows us to generalize the conjecture to the effect that if $R$ is locally an integral domain (that is, $R_M$ is a domain for each maximal ideal $M$), then a nonzero polynomial in $R[X]$ is Gaussian iff its content is locally principal  (see Theorem \ref{main} below).

Our approach is inspired by \cite{GV}. For a finitely generated ideal $I$ of $R$, let $\nu(I)$ be the minimal number of generators of $I$. To prove that $\mathfrak c(f)$ is invertible we first show that $\nu((\mathfrak c f)^n)$ is bounded (Lemma \ref{bounded} below), and conclude that  $\mathfrak c_{R'}(f)=\mathfrak c_{R}R'$ is invertible in $R'$, the integral closure of $R$ (Lemma \ref{int_clos}). To descend from $R'$ to $R$ we simply ``take conjugates'' (see the proof of Theorem \ref{main}).

To bound the number of generators of $(\mathfrak c(f))^n$ we need the following 
proposition (actually, we use just the easier direction $[\Longleftarrow]$).

\begin{prop}\label{power}
Let $f(X)$  be a polynomial in $R[X]$ and let $n\ge1$. Then
$f(X)$ is Gaussian $\iff$ $f(X^n)$ is Gaussian.
\end{prop}

\begin{proof}
Let $g$ by any polynomial in $R[X]$.

$[\implies]$

Write 
$$
g(X)=h_0(X^n)+Xh_1(X^n)+\dots +X^{n-1}h_{n-1}(X^n),
$$ where $h_0(X),\dots,h_{n-1}(X)$ are polynomials in $R[X]$. Since
$f$ is Gaussian, we obtain
\begin{align*}\cn(f(X^n))\cn(g(X))=\cn(f(X^n))\sum_{i=0}^{n-1}\cn(h_i(X^n))=
\sum_{i=0}^{n-1}\cn(f(X^n))\cn(h_i(X^n))\\=
\sum_{i=0}^{n-1}\cn(f(X))\cn(h_i(X))=\sum_{i=0}^{n-1}\cn(f(X)h_i(X))=
\sum_{i=0}^{n-1}\cn(f(X^n)X^ih_i(X^n))\\=
\cn\left(f(X^n)\sum_{i=0}^{n-1}X^ih_i(X^n)\right)
=\cn(f(X^n)g(X)).
\end{align*} Hence $f(X^n)$ is Gaussian.

$[\Longleftarrow]$

Since $f(X^n$) is Gaussian, we obtain
$\cn(f(X)g(X))=\cn(f(X^n)g(X^n))=\cn
(f(X^n))\cn(g(X^n))=\cn(f(X))\cn(g(X))$. Hence $f(X)$ is Gaussian.
\end{proof}

\begin{lemma}\label{bounded}
Let $R$ be a quasilocal domain, let $f(X)$ be a Gaussian polynomial 
in $R[X]$, and let
$I=\cn_R(f)$. Then
$$
\nu(I^n)\le\deg(f)+1
$$
for sufficiently large $n\ge1$.
\end{lemma}

\begin{proof}
By \cite[Corollary 2]{ES} it is enough to show that $\nu(I^{2^m})\le\deg(f)+1$ 
for all $m\ge0$. Let
$$f(X)=g_0(X^2)+Xg_1(X^2),$$ where $g_0(X)$ and $g_1(X)$ are
polynomials in $R[X]$. Since $\cn(f(-X))=\cn(f(X))$ and since $f(X)$ is
Gaussian, we obtain
\begin{align*}
I^2&=(\cn(f))^2=\cn(f(X))\cn(f(-X))=
\cn(f(X)f(-X))\\
&=\cn(g_0(X^2)^2-X^2g_1(X^2)^2)=\cn(g_0(X)^2-Xg_1(X)^2)).
\end{align*}
Since $\deg((g_0(X)^2-Xg_1(X)^2))=\deg(f)$, we infer that
$$\nu(I^2)\le\deg f+1.$$
 Moreover, by Proposition \ref{power}, the polynomial
$g_0(X)^2-Xg_1(X)^2$ is Gaussian since $(g_0(X^2))^2-X^2g_1(X^2)$ is 
a product of two Gaussian polynomials. Thus we may 
proceed by induction on
$m$ to obtain $\nu(I^{2^m})\le\deg f+1$ for all $m\ge0$. This
concludes the proof of the lemma.
\end{proof}

\begin{lemma}\label{int_clos}
Let $R$ be a quasilocal, integral domain, and let $f(X)$ be a Gaussian polynomial
in $R[X]$. Then $\cn_{R'}(f)=R'\cn_R(f)$ is invertible in $R'$.
\end{lemma}

\begin{proof}
By Lemma \ref{bounded}, $\nu(\cn_{R'}(f^n))$ is bounded. Hence, by
\cite[Corollary 1]{ES}, the ideal $\cn_{R'}(f)$ is prestable, and 
by \cite[Lemma F]{ES},  it is an invertible ideal of $R'$ (see also \cite[Theorem 3.1]{GV}).
\end{proof}

\begin{theorem}\label{main}
Let $R$ be a ring which is locally a domain. Then a nonzero
polynomial over $R$ is Gaussian iff its content in $R$ is locally 
principal.
\end{theorem}

\begin{proof}
We may assume that $R$ is quasilocal.  Let $f(X)=\sum_{i=0}^n
a_iX^i$ be
a nonzero Gaussian polynomial in $R[X]$, and let $I=\cn_R(f)$.  By the previous lemma, 
$IR'$ is an invertible ideal in $R'$.
Let $1=\sum_{i=0}^nz_ia_i,$ where $z_i\in (R':I)$ for all $i$.  Let
$g(X)=f(X)\sum_{i=0}^n z_{n-i} X^i=\left(\sum_{i=0}^n
a_iX_i\right)\left(\sum_{i=0}^n z_iX^{n-i}\right)$. Thus $g(X)$ is a
polynomial in $R'[X]$, the coefficient of $X^n$ in $g(X)$ is $1$, and
$f(X)$ divides $g(X)$ in $K[X]$, where $K$ is the fraction field of
$R$:
$$ g(X)=\sum_{i=0}^{2n}\alpha_iX^i,$$ where $\alpha_n=1$.
For each $i\ne n$, there exists a monic polynomial
$h_i$ in $R[X]$ such that $h_i(\alpha_i)=0$; we may decompose all
polynomials $h_i(X)$ into linear factors over some integral extension
$D$ of $R$ containing $R'$:
$$
h_i(X)=\prod_{j=1}^{m_i}(X-\beta_{ij}).
$$ 
Let $\varphi(X)$ be the product of all possible polynomials $\sum_{i=0}^{2n}
\beta_{ij_i}X^i$, where $0\le j_i\le m_i$ for $i\ne n$, and
$j_n=0,\beta_{n,0}=1$. The coefficients of the polynomial $\varphi(X)$
can be expressed as polynomials in the elements $\beta_{ij}$ that are
symmetric in each sequence of indeterminates $X_{i1},\dots,X_{im_i}$
for $i\ne n$. Thus all the coefficients of $\varphi(X)$ are in $R$.
Moreover, $\varphi$ is a product of polynomials in $D[X]$ with unit
content in $D$. Since $D$ is integral over $R$, the polynomial
$\varphi$ has unit content also in $R$. We have $\varphi=f\psi$ for
some polynomial $\psi$ over  $K$. Since the
polynomial $f$ is Gaussian over $R$, we obtain
$R=\cn_R(\varphi)=\cn_R(f)\cn_R(\psi)$, thus $\cn_R(f)$ is an invertible ideal in $R$.
\end{proof}

Theorem \ref{main} implies that the Gaussian property of a polynomial
over an integral domain depends just on its content.  In the next
corollary we present further immediate consequences of Theorem
\ref{main}.

\begin{cor}\label{conseq}
Let $R$ be an integral domain. We have
\begin{enumerate}
\item
If $f$ and $g$ are polynomials in $R[X]$ with the same content, then
$f$ is Gaussian iff $g$ is Gaussian. In particular, a polynomial
obtained from a Gaussian polynomial over $R$ by permuting its
coefficients is Gaussian.
\item
A Gaussian polynomial over $R$ is Gaussian over any ring
extension of $R$.
\item
All  polynomials over $R$ are Gaussian iff $R$ is a Pr\"ufer domain (this result was already obtained in a more general form in \cite[Theorem 1.3]{AK}).
\end{enumerate}
\end{cor}

\begin{example}[cf. Corollary \ref{conseq} (2)]
We consider the ring $R=k[s,t]/(s,t)^2$, where $k$ is a field, and $s$
and $t$ are independent indeterminates over $k$, thus all polynomials over $R$ 
are Gaussian (cf. \cite{{AK},{GV}}).
However, if $u$ and $v$ are indeterminates
over $R$, then the polynomial $s+tX$ is not Gaussian over the ring extension
$R[u,v]$. Indeed, $sv\in\cn_{R[u,v]}(s+tX)\cn_{R[u,v]}(u+vX)$, but
$st\notin \cn_{R[u,v]}((s+tX)(u+vX))$, that is, $sv\notin (su, tv,
sv+tu)$.  \qed
\end{example}

\begin{example}
 A factor of a nonzero Gaussian polynomial over an integral domain
is not necessarily Gaussian. Moreover, if $R$ is an integral domain,
$f(X)\in R[X]$, and $f^2$ is Gaussian, then $f$ is not necessarily
Gaussian.

Indeed, let $k$ be a field of characteristic $2$, and let
$R=k[s,t,\frac{s^2}{t^2}]$, where $s$ and $t$ are indeterminates over
$k$. Let $f(X)=s+tX$. We have $f^2=s^2+t^2X^{2}$ is Gaussian since its
content is generated by one element, namely, by $t^2$, but
$f$ is not Gaussian since $st\in(\cn (f))^2\smallsetminus \cn(f^2)$.
\qed
\end{example}
However, a factor of a nonzero Gaussian polynomial over an integrally closed domain is 
Gaussian (a domain $R$ is integrally closed iff the multiplicative set of the polynomials with invertible content, that is, the set of nonzero Gaussian polynomials, is saturated).

We conjecture that Theorem \ref{main} can be extended to reduced rings,
to the effect that a polynomial over a reduced ring is Gaussian iff its content is locally principal. However, if $R=k[s,t]/st$, where $k$ is a field, then the polynomial 
$f(X)=\bar s+\bar tX\in R[X]$ is not Gaussian, although $R$ is a principal ideal domain modulo each of its two nonzero minimal primes, namely, $(\bar s)$ and $(\bar t)$. To show that the polynomial $f(X)$ is not Gaussian, let $g(X)= \bar t+\bar sX$, thus
$fg=(\bar s^2+\bar t^2)X$; we have $\bar s^2\in\cn(f)\cn(g)$, but
$\bar s^2\notin\cn(fg)$ since $s^2\notin (st, s^2+t^2)k[s,t]$ (cf. \cite{GV}).

Finally, we conjecture that Theorem \ref{main} can be generalized to any number of indeterminates over any reduced ring $R$; it is enough to consider the case of a finite number of indeterminates. By the above proof, the conjecture holds for an integral domain $R$ of finite characteristic; more generally, it is enough to assume that the residue fields $R/M$ for $M$ a maximal ideal in $R$, are of finite characteristic (see \cite{GV}).

\end{document}